\documentclass[12pt,fleqn,legno]{article}
\usepackage{amsmath}
\usepackage{amsfonts}
\usepackage[T1]{fontenc}
\usepackage{amsthm}
\usepackage{mathrsfs}

\newcommand{\Rset}{\mathbb{R}}
\newcommand{\Cset}{\mathbb{C}}

\newtheorem{thm}{THEOREM}

\newtheorem{cor}{COROLLARY}
\newtheorem{rem}{Remark}

\begin{document}
$10^{th}$ Vilnius Conference on Probability

and Mathematical Statistics

Vilnius, June 28-July 2, 2010.

\medskip
\medskip
Zbigniew J. JUREK

\medskip
\begin{center}
{\large \textbf{The Random Integral Representation Conjecture}:  \ \
\qquad \qquad a quarter of a century later}
\end{center}

\begin{center}
\end{center}

\begin{quote}
\textbf{Abstract.} In Jurek 1985 and 1988 the random integral
representations conjecture was stated. It claims that (some) limit
laws can be written as probability distributions of random integrals
of the form $\int_{(a,b]}h(t)dY_{\nu}(r(t))$, for some deterministic
functions $h$, $r$ and a L\'evy process $Y_{\nu}(t),t\ge 0$. Here we
review situations where a such claim holds true. Each theorem is
followed by a remark which gives references to other related papers,
results as well as some historical comments. Moreover, some open
questions are stated.



\emph{Abbrivated title:} Integral Representation Conjecture
\end{quote}

\newpage
In 1985 in  \emph{The Annals of Probability, vol. 13, No.
2, on the page 607}  and later on in 1988 in \emph{Probability
Theory and Related Fields, vol. 78, on the page 474},  it was
conjectured that:

\medskip
\medskip
\textbf{\emph{Each class of limit distributions, derived from
sequences of independent random variables, is the image of some
subset of ID by some mapping defined as a random integral.}}
\footnote{Chatterji's subsequence principle claiming that: \
\emph{Given a limit theorem for independent identically distributed
random variables under certain moment conditions, there exists an
analogous theorem such that an arbitrary-dependent sequence (under
the same moment conditions) always contains a subsequence satisfying
this analogous theorem}, was proved by David J. Aldous (1977).
Although we do not expect that  the above Conjecture and Chatterji's
subsequence principle are mathematically related, however, one may
see a $"$philosophical$"$ relation between those two.}

\medskip
\medskip
\noindent More formally, one claims that for a class $K$ of limiting
probability distributions on a Banach space $E$ there exist: a
function $h$ (a space scaling), a function $r$ (a time change), an
interval $A$ (in a positive half-line) and a subset $\mathcal{D}$ of
$ID$ (the class of all infinitely divisible distributions) such that
\begin{equation}
K\equiv I^{h, r}_A(\mathcal{D}):=\{I^{h,
r}_A(\nu):=\mathcal{L}\big(\int_A \, h(s)\,dY_{\nu}(r(s))\big): \nu
\in \mathcal{D}\},
\end{equation}
where $Y_{\nu}(s), s\ge 0,$ is an E-valued L\'evy process with
cadlag paths such that its probability distribution at time 1,
$\mathcal{L}(Y_{\nu}(1))=\nu$ and $\mathcal{D}$ denotes the domain
of existence of the above random integral; cf.
www.math.uni.wroc.pl/$\sim$zjjurek  \, \,( The Conjecture.)

Note that the notation $I^{h, r}_A$  for the random integral
transformation can be simplified as follows
\[
I^{h,\, r}_A\,(\nu)=\int_0^{\infty}\,1_A\,(r^*(s))
h(r^*(s))\,dY_{\nu}(s) =
I^{\,\tilde{h}(s),\,\,s\,}_{(0,\infty)}(\nu) \equiv
I^{\tilde{h}}(\nu), \ \ \ \qquad (1^{\prime})
\]
where $\tilde{h}(s):=1_A\,(r^*(s)) h(r^*(s))$ and $r^*$ is the
inverse function of $r$.

\medskip
The term \emph{random integral} emphasizes  the fact that the
integrand $h$ is a deterministic function. Thus for $A=(a,b]$ we may
define the random integral by the formal integration by parts
formula, i. e.,
\begin{equation}
\int_A h(s)\,dY_{\nu}(r(s)):=h(b)Y_{\nu}(r(b))-
h(a)Y_{\nu}(r(a))-\int_A\,Y_{\nu}(r(s))dh(s), \ \
\end{equation}
provided $h$ is of a bounded variation. Thus approximating the
right-hand side integral, by Riemann -Stieltjes sums,  we get the
formula for the Fourier transform
\begin{equation}
\log\widehat{(I^{h,r}_A(\nu))}(y)=\int_A\log\widehat{\nu}(h(s)y)dr(s)
, \ y \in E^{'} \mbox{is the dual Banach space.}
\end{equation}
Random integrals on half-lines $(a,\infty)$ are defined as weak
limits of the integrals (1) for $(a,b]$ as $b\to \infty$.

Below we review the old results as well as the more recent ones. In
remarks after  each theorem we point out to other related facts and
papers. This survey is divided into three basic parts.  While the
last one rephrases the conjecture.

\medskip
\medskip
\textbf{1. From a class of limit laws to a class of random
integrals.}

\medskip
\textbf{(a)} \ For the L\'evy class $L$ of selfdecomposable
probability measures, that coincides with the class of limiting
distributions of the following of sequences
\begin{equation}
T_{a_{n}}(\xi_1+\xi_2+...+\xi_n) + x_n,  \ \ \ T_a(x):=a\,x,  \ \
a>0, \ x\in E,
\end{equation}
where $(\xi_i)$ are independent E-valued random variables, $x_n\in
E$, $a_n>0$ and the summands in (4) are uniformly infinitesimal we
have:
\begin{thm} (Jurek and Vervaat(1983)).
For the class $L$ we have that
\[
L=\{\, I^{\,e^{-s},\, s}_{(0,\infty)}\,(\nu):\ \nu\in
ID_{\log}\,\}=\{\mathcal{L}(\int_{(0,\infty)}e^{-s}dY_{\nu}(s)):\nu\in
ID_{\log}\},
\]
where $ID_{\log}$ is the class of all infinitely divisible measures
on $E$ that integrate the function $\log(1+||x||)$.
\end{thm}

\begin{rem}
\emph{ (i) \ S. J. Wolfe (1982) and K. Sato with M. Yamazato (1984)
had similar characterizations but with proofs valid \underline{only}
in Euclidean spaces. (ii) To the processes $Y_{\nu}$ above and more
generally to ones in (1) we refer to as \emph{ the background
driving L\'evy processes; in short: BDLP}; cf. Jurek (1996). (iii) A
connection between selfdecomposable distributions and the one
-dimensional Ising models in statistical physics were shown in Jurek
(2001). (iv) Replacing $T_a$'s in (4) by \underline{arbitrary linear
operators} we get so called \emph{operator-limit distributions
theory}; cf. Jurek and Mason (1993) or Meerschaert and Scheffler
(2001). See also Urbanik (1972, 1978)}.
\end{rem}

 \textbf{(b)} \ If one assumes that $L_1 \equiv L$ and for positive integer $m\ge2$ one
 defines
the class $L_m$ as a class of limits of (4) but such that
$\mathcal{L}(\xi_i)\in L_{m-1}$ then $L_{m+1}\subset L_m$ and
moreover we have :
\begin{thm} (Jurek(1983)).
For the class $L_m , m=1,2,...$ we have that
\[
L_m=\{\, I^{\,e^{-s},\, s^m/m!}_{(0,\infty)}\,(\nu):\ \nu\in
ID_{\log^m}\,\}=\{\mathcal{L}(\int_{(0,\infty)}e^{-s}dY_{\nu}(\frac{s^m}{m!}):\nu\in
ID_{\log^m}\},
\]
where $ID_{\log^m}$ is the class of all infinitely divisible
measures on $E$ that integrate the function $\log^m(1+||x||)$.
\end{thm}

\begin{rem}
\emph{(i) The idea of classes $L_m$ belongs to Urbanik (1973) with a
different scheme of summation; see also Kumar and Schreiber (1979).
The iterative approach was proposed in K. Sato (1980) (for Euclidean
spaces) and later generalized in Jurek (1983a) in two direction:
replacing Euclidean space by an arbitrary separable Banach space $E$
and \underline{more importantly} replacing the group $(T_a, a>0)$ of
dilation by  an arbitrary strongly continuous one-parameter group
$\mathbb{U}$ of bounded linear operators on $E$. \noindent(ii) The
particular case of the group $\mathbb{U}:=\{e^{-tQ}: t\in \Rset\}$,
where $Q$ is a fixed bounded linear operator on a Banach space $E$,
was investigated in Jurek (1983b), where  it was shown that
$L_m(Q)=\{\,I^{e^{-tQ}, \, t}_{(0,\infty)}\,(\nu): \nu \in
ID_{\log^m} \}$, i.e., in Theorem 2 the scalar function $e^{-t}$ is
replaced by the operator-valued function $e^{-tQ}$. Here one needs a
new norm and the polar coordinates in a Banach space; cf. Jurek
(1984). (iii) N. Thu (1986) extended classes $L_m, \ m=1, 2,...,$ to
$L_{\alpha}, \alpha>0,$ by using the fractional calculus. }
\end{rem}

\medskip
\textbf{(c)} \ Let us replace the linear normalization in (4)  by
the \underline{\emph{non-linear}} shrinking s-operation $U_r$ (r>0)
and consider the class $\mathcal{U}$ of limiting distributions in
the following scheme
\begin{equation}
U_{r_{n}}(\xi_1)+ U_{r_{n}}(\xi_2)+...+U_{r_{n}}(\xi_n) + x_n,  \
U_r(x):= \max(||x||-r,0)\frac{x}{||x||}, x\neq 0,
\end{equation}
where the summands are uniformly infinitesimal. Limiting
distributions of (5) are called \emph{s-selfdecomposable
distributions.}
\begin{thm} (Jurek(1985)).
For the class $\mathcal{U}$ of s-selfdecomposable distributions we
have that
\[
\mathcal{U}=\{\, I^{\,s,\, s}_{(0,1)}\,(\nu):\ \nu\in ID
\,\}=\{\mathcal{L}(\int_{(0,1)}s\,dY_{\nu}(s)):\nu\in ID\},
\]
where $ID$ is the class of all infinitely divisible measures .
\end{thm}

\begin{rem}
\emph{(i) Note that $U_r(U_s(x))=U_{r+s}(x)$ (semigroup of
non-linear transformations) and for positive random variable
$\xi>0$, we get $U_r(\xi)=(\xi-r)^+$ (the positive part), i.e., it
coincides with the famous financial derivative \emph{call option}.
(ii) Characterizations of the s-selfdecomposable distributions in
terms of the L\'evy-Khintchine formula were presented during
\underline{$2^{nd}$ Vilnius Conference}; cf. Jurek (1977). Complete
proofs were given in Jurek (1981). (iii) The CLT for s-operations
$U_r$ was proved by Housworth and Shao (2000). (iv) Classes
 $\mathcal{U}_{\beta}: = I^{\,s, \,s^{\beta}}_{(0,1)}(ID)$, were
investigated in a series of papers: Jurek (1988) for  $\beta>0$,
Jurek (1989) for $-1\le \beta<0$ and Jurek and Schreiber (1992) for
$-2<\beta\le -1$. In the last two cases, the stable distributions
appeared as convolution factors of the limiting distributions.
Measures from $\mathcal{U}_{\beta}$ are called \emph{generalized
s-selfdecomposable distributions.}}
\end{rem}

In a similar way as the classes $L_m$ were introduced in Theorem 2,
one may iterate  the random integral mapping $I^{\,s,\, s}_{(0,1)}$
from Theorem 3, and get the classes $\mathcal{U}^{<m>}$ for which we
have
\begin{thm} (Jurek(2004)).
For the class $\mathcal{U}^{<m>}$  ( with $m=1,2, ... $) of  m-times
s-selfdecomposable distributions we have that
\begin{multline*}
\mathcal{U}^{<m>}=\{\, I^{\,s,\, \tau_m(s)}_{(0,1)}\,(\nu):\ \nu\in
ID \,\}=\{\mathcal{L}(\int_{(0,1)}s\,dY_{\nu}(\tau_m(s)):\nu\in
ID\}, \\
\tau_m(s): =\frac{1}{(m-1)!}\int_0^s (-\log u)^{m-1}du,  \ \ \  \ \
0<s \le 1 ,  \ \ \ \ \ \ \ \ \ \
\end{multline*}
where $ID$ is the class of all infinitely divisible measures .
\end{thm}

Although the classes $L_m$ and $\mathcal{U}^{<m>}$ originated in two
different limiting schemes ( via the linear dilations $T_a$  and the
non-linear s-operations $U_r$, respectively) they  still admit some
unexpected relations.

\begin{cor} (Jurek (2004))
(a) We have inclusions

 \ \ \ \ \ \ $L_{m+1}\subset  \mathcal{U}^{<m+1>}\subset \mathcal{U}^{<m>}\subset
ID , \ \ \ \ \ m=1,2,...$.

(b)  \ \ \ $L_{\infty}:= \bigcap_{m=1}^{\infty}\, L_m=
\mathcal{U}^{<\infty>}:=\bigcap_{m=1}^{\infty}\mathcal{U}^{<m>} \\
 = \mbox{the smallest closed convolution semigroup that contains all stable
measures.}$
\end{cor}

\begin{rem}
\emph{(i)  \ Maejima and Sato (2009) proved that besides the two
instances described in Corollary 1 (b)  there are another three
classes for which infinitely many integral iterations  lead to the
smallest closed convolution semigroup  that contains  all stable
measures. (ii) \underline{ Still an open question} is to describe
$L_{\infty}(\mathbb{U}):=\bigcap_{m=1}^{\infty} L_m(\mathbb{U})$,
where $\mathbb{U}$ is an one-parameter group of bounded linear
operators on a Banach space $E$; comp. Remark 2(i) and Jurek
(1983).}
\end{rem}

\medskip
\medskip
\textbf{2. From a class of random integrals to ... }

\medskip
The original aim (in the 80's of the last century) was to identify
\emph{a given class K of limit distributions} as a collection of
probability distributions of some random integrals; comp. the above
Section 1. Later on, more often questions were asked whether
\emph{given class of distributions} (or Fourier transforms or L\'evy
spectral measures), can be described in terms of some random
integrals. In this section we discuss only two of  such examples.

\medskip
\textbf{(a)} \ D. Voiculescu and others studying so called
\emph{free-probability} introduced new binary operations on
probability measures and termed them \emph{free-convolutions}; cf.
Bercovici-Voiculescu (1993) and references therein. For the
\emph{additive free-convolution} $\Box$ the Voiculescu transform
$V_{\nu}(z), z \in \Cset,$ (an analogue  of the characteristic
function $\hat{\nu}(t), t\in\Rset$) is additive. Namely, $V_{\nu_1
\Box \nu_2}(z)=V_{\nu_1}(z)+V_{\nu_2}(z)$. This property allowed to
introduce a notion of free-infinite divisibility.

\begin{thm} (Jurek (2007).)
A probability measure $\nu$ is $\Box$-infinitely divisible if and
only if there exist a unique $\ast$-infinitely divisible probability
measure $\mu$ such that
\begin{equation*}
(it)\,V_{\nu}((it)^{-1})=\log\big( I^{s,\,
1-e^{-s}}_{(0,\infty)}(\mu)\big)^{\widehat{}}(t)= \log
\Big(\mathcal{L}
(\int_0^{\infty}sdY_{\mu}(1-e^{-s}))\Big)^{\widehat{}}(t), \ t\neq
0,
\end{equation*}
where $(Y_{\mu}(t), t\ge 0)$ is a L\'evy process such that
$\mathcal{L}(Y_{\mu}(1))=\mu$.
\end{thm}

\begin{rem}
\emph{Using Theorem 5 one can easily see that for we have the
integral mapping
\begin{equation}
\mathcal{K}(\mu):=I^{s, \, 1-e^{-s}}_{(0,\infty)}(\mu)=I^{\,-\log s,
\, s}_{(0,1)}(\mu), \ \ \ \mu \in ID.
\end{equation}
Mapping $\mathcal{K}$ was called the $\Upsilon$ (upsilon) transform
and studied from  a different point of view by Barndorff-Nielsen ,
Maejima and Sato (2006); also by  Barndorff-Nielsen, Rosiñski and
Thorbjorsen (2008), and by Maejima and Sato (2009).}
\end{rem}

\medskip
\textbf{(b) }\ \ Thorin class $T(\mathbb{R}^d)$ is an example of the
class of infinitely divisible distributions defined by properties of
their L\'evy spectral measures that, later on, was characterized by
some random integrals;  for more details cf. Maejima and Sato
(2009), p. 121; for related results  cf. Grigelionis (2007). For the
class $G$ distributions see Aoyama and Maejima (2007).

\begin{thm}(Maejima and Sato (2009))
\begin{equation*}
T(\mathbb{R}^d)=\{\mathcal{L}(\int_0^{\infty}e^*(t)dY_{\mu}(t)):
\mu\in ID_{\log} \}=\{ I^{\,s,\, e(s)}_{(0,\infty)}(\mu):\mu\in
ID_{\log}\},
\end{equation*}
where $e(s):=\int_s^{\infty}u^{-1}e^{-u}du, s>0$ and $e^*(t), t>0$
is its inverse function.
\end{thm}
\begin{rem}
\emph{In Jurek (2007), Proposition 4, the class
$\mathcal{TS}_{\alpha}$ of tempered stable distributions with the
index $0<\alpha<1$ was identified as the class of random integrals
\[
\mathcal{TS}_{\alpha}=\{ I^{\,s,\,\Gamma(-\alpha,
s)}_{(0,\infty)}(\mu):\mu\in ID_{\alpha}\}, \ \ \mbox{and} \ \
\Gamma(-\alpha, s):=\int_s^{\infty}\,w^{-\alpha-1} e^{-w}dw ,\ \
s>0,
\]
and $ID_{\alpha}$ denotes the set of all infinitely measures whose
L\'evy spectral measures integrate $||x||^{\alpha}$ over the space
$E$.}
\end{rem}

\medskip
\medskip
\textbf{3. A calculus on (L\'evy exponents of) ID distributions. }

\medskip
On the random integrals (1), (viewed as mappings defined on (some)
infinitely divisible probability measures $\mu$), one can perform
transformations such as compositions or the arithmetic operations.
Because of the formula (3) all of that operations have natural
generalizations to the L\'evy exponents, that is, the logarithms
$\Phi:=\log\hat{\mu}$ of Fourier transforms of $ID$ measures $\mu$.
Such a calculus may lead to new factorization properties. For the
simplicity of the notations let us pu
\begin{equation}
\mathcal{I}(\mu)\equiv I^{\,e^{-s}\,, s}_{(0,\infty)}(\mu), \
\mbox{for} \ \mu\in ID_{\log} \ \ \ \mbox{and} \   \
\mathcal{J}(\mu)\equiv I^{s\,, s}_{(0, 1)}(\mu), \ \mbox{for} \ \
\mu\in ID.
\end{equation}
Then we have
\begin{thm} (Jurek (1985) and (2008))\ (i) For  the mappings $\mathcal{I}$ and $\mathcal{J}$
and $\nu\in ID_{\log}$ we have the identity
\begin{equation*}
\mathcal{I}(\nu\ast\mathcal{J}(\nu))=\mathcal{I}(\nu)\ast\mathcal{I}(\mathcal{J}(\nu))=
\mathcal{I}(\nu)\ast\mathcal{J}(\mathcal{I}(\nu))=
\mathcal{J}(\nu);
\end{equation*}
(ii) For each selfdecomposable measure $\mu\in L$ there exists a
unique \\ s-selfdecomposable measure  $\tilde{\mu}\in\mathcal{U}$
such that
\begin{equation}
\mu=\tilde{\mu}\ast\mathcal{I}(\tilde{\mu})  \ \ \ \mbox{and} \ \
\mathcal{J}(\mu)=\mathcal{I}(\tilde{\mu}).
\end{equation}
\end{thm}

\begin{rem}
\emph{(i) When one considers $\mathcal{I}$ and $\mathcal{J}$ as the
mappings defined on the L\'evy exponents $\Phi$'s (via the equation
(3)) then one gets
$\mathcal{I}(\mathcal{J})=\mathcal{I}-\mathcal{J}$ or
$\mathcal{J}(I+\mathcal{I})=\mathcal{I}$ or
$\mathcal{I}(I-\mathcal{J})=\mathcal{J}$ or
$(I-\mathcal{J})(I+\mathcal{I})=I$. (ii) Iterating the property (8)
we get a convolution decomposition of selfdecomposable distributions
with m-times s-selfdecomposable distributions (from the class
$\mathcal{U}^{<m>}$) as the factors; cf. Jurek (2008).}
\end{rem}

From Theorem 1 and formula (7) we have that
$L=\mathcal{I}(ID_{\log})$. We say that a selfdecomposable
$\mu=\mathcal{I}(\rho)$ has \emph{the factorization property} if
$\mathcal{I}(\rho) \ast\rho\in L$ and let $L^f$ denotes the class of
all class $L$ distributions having the factorization property.
\begin{thm} (Ikasnov, Jurek and Schreiber (2004);
Czy¿ewska-Jankowska and Jurek (2009)) \ (i) \
$L^f=\mathcal{I}(\mathcal{J}(ID_{\log}))$;
\begin{equation*} (ii) \
L^f =
I^{e^{-s},\,s+e^{-s}-1}_{(0,\infty)}(ID_{\log})=\big\{\mathcal{L}\bigl(\int_{0}^{\infty}
e^{-s}\,dY_{\nu}(s +e^{-s}-1))\,: \ \nu\in
 ID_{\log}\big\}.
\end{equation*}
\end{thm}

\begin{rem}
\emph{(i) Note that $ I^{e^{-s},\,s+e^{-s}-1}_{(0,\infty)}= I^{\, s
,\,s-\log s-1}_{(0,1)}$.  (ii) One of the most important examples of
the class $L^f$ distribution is the L\'evy's stochastic area
integral (the hyperbolic sine  characteristic function). In Jurek
and Yor (2004) BDLP were identified as Bessel squared processes.
(iii) Above we have that $I^{s,s}_{(0,1)}(I^{e^{-s},
s}_{(0,\infty)})=I^{e^{-s},\,s+e^{-s}-1}_{(0,\infty)}$, it means
that a composition of two random integrals is again a random
integral.}
\end{rem}

\medskip
\medskip
\textbf{4. Concluding remarks.}

Taking into account the above historical survey :

$(\alpha)$ \ \  Can we still hope for a general proof of the RANDOM
INTEGRAL REPRESENTATION CONJECTURE ?

\medskip
Even without settling the previous question:

$(\beta)$ \ \  Can we develop $"$an abstract theory$"$ of a calculus
on random integral mappings $I^{h,r}_A$ (or on the corresponding
L\'evy exponents or L\'evy spectral measures)?

\medskip
Since (many) classes, of probability distribution of random integral
mappings discussed above, naturally form convolution semigroups:

$(\gamma)$ \ \ Can we find  structural descriptions of ALL (closed)
convolution subsemigroups of the semigroup $ID$ ( of all infinitely
divisible distributions)?

\medskip
\begin{center}
\textbf{References}
\end{center}
\noindent [1] D. J. Aldous (1977), Limit theorems for subsequences
of arbitrary-dependent sequences of random variables, \emph{Z.
Wahrscheinlichkeitstheorie verw. Gebiete}, vol. 40, pp.59-82.

\noindent[2] T. Aoyama and M. Maejima (2007), Characterizations of
subclasses of type G distributions on $\Rset^d$ by stochastic random
integral representation, \emph{Bernoulli}, vol. 13, pp. 148-160.

\noindent[3]  O. E. Barndorff-Nielsen, M. Maejima and K. Sato
(2006), Some classes of multivariate infinitely divisible
distributions admitting stochastic integral representations,
\emph{Bernoulli}, \textbf{12}, 1-33.

\noindent[4]  O. E. Barndorff-Nielsen, J. Rosiñski and S.
Thorbjornsen (2008), General $\Upsilon$-transformations,
\emph{ALEA}, \textbf{4}, 131-165.

 \noindent[5]  H. Bercovici and D. V. Voiculescu (1993), Free
convolution of measures with unbounded support, \textit{Indiana
Math. J.} \textbf{42}, 733-773.

\noindent[6]  A. Czy¿ewska-Jankowska and Z. J. Jurek (2008),
Factorization property of s-selfdecomposable measures and class
$L^f$ distributions, \textit{to appear}.

\noindent[7] B. Grigelionis (2007), Extended Thorin classes and
stochastic integrals, \emph{Liet. Matem. Rink.} \textbf{47} ,
497-503.

\noindent[8]  E. Housworth and Q-M. Shao (2000), On central limit
theorems for shrunken random variables, \emph{ Proc. Amer. Math.
Soc.} \textbf{128}, No 1, 261-267.

\noindent[9]  A. M. Iksanov, Z. J. Jurek and B.M. Schreiber (2004),
A new factorization property of the selfdecomposable probability
measures, \emph{Ann. Probab.} \textbf{32}, No 2, 1356-1369.

\noindent[10] Z. J. Jurek (1977), Limit distribution for sums of
shrunken random variables. In: \emph{Second Vilnius Conf. Probab.
Theor. Math. Stat.}, Abstract of Communications \textbf{3}, 95-96.

\noindent[11]  Z. J. Jurek (1981), Limit distributions for sums of
shrunken random variables, \textit{Dissertationes Mathematicae},
vol. \textbf{185}, PWN Warszawa.

\noindent[12]  Z. J. Jurek (1983a) , Limit distributions and
one-parameter groups of linear operators on Banach spaces,
\textit{J. Multivariate Anal.} \textbf{13}, 578-604.

\noindent[13]  Z. J. Jurek (1983b), The classes $L_{m}(Q)$ of
probability measures on Banach spaces, \textit{Bull. Acad. Polon.
Sci., S\'er. Sci. Math.} \textbf{31}, 51-62.

\noindent[14] Z. J. Jurek (1984), Polar coordinates in Banach
spaces, \textit{Bull. Polish Acad. Sci; Math.} \textbf{32}, 61- 66.

\noindent[15]  Z. J. Jurek (1985), Relations between the
$s$-selfdecomposable and selfdecomposable measures, \textit{Ann.
Probab.} \textbf{13}, 592-608.

\noindent[16]  Z. J. Jurek  (1988). Random Integral representation
for Classes of Limit Distributions Similar to L\'evy Class $L_{0}$,
 \emph{Probab. Th. Fields.} \textbf{78},  473-490.

\noindent[17]  Z. J. Jurek  (1989). Random Integral representation
for Classes of Limit Distributions Similar to L\'vy Class $L_{0}$,
II, \emph{Nagoya Math. J.}, \textbf{114}, 53-64.

\noindent[18]  Z. J. Jurek (1996), Series of independent exponential
random variables. In \textit{Proc.7th Japan-Russia Symposium on
Probab. Theor. Math. Statist., Tokyo 1995}, S. Watanabe, M.
Fukushima, Yu. V. Prokhorov and A. N. Shiryaev, eds., World
Scientific, Singapore, 174-182.

\noindent[19] Z. J. Jurek (2001), 1D ising models. compund geometric
distributions and selfdeomposability, \emph{Reports on Math.
Physics}, \textbf{47}, no 1, 21-30.

\noindent[20] Z. J. Jurek (2004), The random integral representation
hypothesis revisited: new classes of s-selfdecomposable laws. In:
ABSTRACT AND APPLIED ANALYSIS, \emph{Proc. International Conf.
Hanoi, Vietnam, 13-17 August 2002}, pp. 495-414, World Scientific,
Hongkong (2004).

\noindent[21]  Z. J. Jurek (2007), Random integral representations
for free-infinitely divisible and tempered stable distributions,
\emph{Stat.$\&$ Probab. Letters} \textbf{77}, 417-425.

\noindent[22]  Z. J. Jurek (2008), A calculus on L\'evy exponents
and selfdecomposability on Banach spaces, \emph{Probab. Math. Stat.}
\textbf{28}, Fasc. 2, 271-280.

\noindent[23]  Z. J. Jurek (2009), On relations between Urbanik and
Mehler semigroups, \emph{Probab. Math. Stat.} \textbf{29}, Fasc. 2,
297-308.

\noindent[24]  Z. J. Jurek and J. D. Mason (1993),
\emph{Operator-limit distributions in probability theory}, J. Wiley
and Sons, New York.

\noindent[25] Z. J. Jurek and B. M. Schreiber (1992), Fourier
transforms of measures from classes $\mathcal{U}_{\beta}$, \
$-2<\beta <-1$, \emph{J. Multivar. Analysis}, \textbf{41}, 194-211.

 \noindent[26]  Z. J. Jurek and W. Vervaat (1983), An integral
representation for selfdecomposable Banach space valued random
variables, \textit{Z. Wahrsch. verw. Gebiete} \textbf{62}, 247-262.

\noindent[27] Z. J. Jurek and M. Yor (2004), Selfdecomposable laws
associated with hyperbolic functions, \emph{Probab. Math. Stat.}
\textbf{24}, no.1, pp. 180-190.

\noindent[28] A. Kumar and B. M. Schreiber (1979), Representation of
certain infinitely divisible probability measures on Banach spaces,
\textit{J. Multivariate Anal.} \textbf{9}, 288--303.

\noindent[29] M. Maejima and K. Sato (2009), The limits of nested
subclasses of several classes of infinitely divisible distributions
are identical with the closure of the class of stable distributions,
\emph{Probab. Th. Rel. Fields,} \textbf{145}, 119-142.

\noindent[30] M. M. Meerschaert and H. P. Scheffler (2001),
\emph{Limit distributions of sums of independent random vectors}, J.
Wiley and Sons, New York.

\noindent[31] K. Sato (1980), Class $L$ of multivariate
distributions and its subclasses, \textit{J. Multivariate Anal.}
\textbf{10}, 201-232.

\noindent[32] K. Sato and M. Yamazato (1984),
Operator-selfdecomposable distributions as limit distributions of
processes of Ornstein-Uhlenbeck type, \emph{Stoch. Processes and
Appl.} \textbf{17}, 73-100.

\noindent[33]  N. van Thu (1986), An alternative approach to
multiply selfdecomposable probability measures on Banach spaces,
\textit{Probab. Theory and Rel. Fields} \textbf{72}, 35-54.

\noindent[34] K. Urbanik (1972), L\'evy's Probability measures on
Euclidean spaces, \emph{Studia Math.} \textbf{44}, 119-148.

\noindent[35] K. Urbanik (1973), Limit laws for sequences of normed
sums satisfying some stability conditions. In \textit{Multivariate
Analysis III}, P.R. Krishnaiah, ed., Academic Press, New York,
225--237.

\noindent[36] K. Urbanik (1978), L\'evy's Probability measures on
Banach spaces, \emph{Studia Math.} \textbf{63}, 283--308.

\noindent[37] S. J. Wolfe (1982), On a continuous analogue of the
stochastic difference equation $X_n=\rho X_{n-1}+B_n$, \emph{Stoch.
Processes and Appl.} \textbf{12}, 301-312.

\medskip
\medskip
\medskip
\noindent Institute of Mathematics, University of Wroclaw \\ Pl.
Grunwaldzki 2/4, 50-384 Wroclaw, Poland \\
e-mail: zjjurek@math.uni.wroc.pl\ \ \
www.math.uni.wroc.pl/$\sim$zjjurek

\end{document}